\documentclass[letterpaper,10pt,conference]{ieeeconf}

\IEEEoverridecommandlockouts % This command is only needed if you want to use the \thanks command
\usepackage{cite}
\usepackage{graphicx}
\usepackage{color}
\usepackage{amsmath,amssymb,bm}
\usepackage{tikz-cd}
\usepackage{hyperref}
\hypersetup{
    colorlinks=true,
    linkcolor=blue,
    filecolor=magenta,      
    urlcolor=cyan,
}
\usepackage{algorithm, caption}
\usepackage[noend]{algpseudocode}
\usepackage{placeins}
\usepackage{subfig}
\usepackage{wrapfig}

\newcommand{\R}{\mathbb{R}}
\newcommand{\E}{\mathbb{E}}
\newcommand{\C}{\mathbb{C}}

\newcommand{\koop}{\mathcal{K}}
\DeclareMathOperator{\Tr}{Tr}

\newcommand{\T}{\textsf{T}}

\title{\LARGE \bf
An MCMC Method for Uncertainty Set Generation\\ via Operator-Theoretic Metrics
}
\author{Anand Srinivasan and Naoya Takeishi% <-this % stops a space
\thanks{%
A. Srinivasan is with the Mathematics Department at the Massachusetts Institute of Technology, Cambridge, MA, USA
{\tt\small asrini@mit.edu}%
}%
\thanks{N. Takeishi is with RIKEN Center for Advanced Intelligence Project, Tokyo, Japan
{\tt\small naoya.takeishi@riken.jp}%
}%
}

\begin{document}

\maketitle
\thispagestyle{empty}
\pagestyle{empty}

\begin{abstract}
Model uncertainty sets are required in many robust optimization problems, such as robust control and prediction with uncertainty, but there is no definite methodology to generate uncertainty sets for nonlinear dynamical systems. In this paper, we propose a method for model uncertainty set generation via Markov chain Monte Carlo. The proposed method samples from distributions over dynamical systems via metrics over transfer operators and is applicable to general nonlinear systems. We adapt Hamiltonian Monte Carlo for sampling high-dimensional transfer operators in a computationally efficient manner. We present numerical examples to validate the proposed method for uncertainty set generation.
\end{abstract}

\newtheorem{assumption}{Assumption}

%==================================================

\section{INTRODUCTION}

% Background: Model uncertainty set generation.
Generating \emph{model uncertainty sets} of dynamical systems is a universal problem in situations such as robust control, prediction with uncertainty, and scenario optimization.
For example, in a min-max model predictive control (MPC) problem, one would like to solve 
\begin{equation*}
    \min_{\{u\}} \max_{F \in \mathcal{F}} ~ J(\{u\}, x_0, F),
\end{equation*}
where $J$ is some loss function, $u$ denotes control signals, $x_0$ is an initial state, and $F$ is a dynamics model.
Here, $\mathcal{F}$ is a model uncertainty set with which the worst-case performance is to be optimized and thus is key to good performance of the robust controller.
However, the way to configure a good $\mathcal{F}$ is not trivial in general.

% Existing method: Bayesian inference.
A possible strategy for model uncertainty set generation is via Bayesian inference, with which we can create an uncertainty set from an inferred posterior.
One of the common difficulties in Bayesian inference is that we need to prepare appropriate priors and observation models, which is sometimes challenging when the target dynamics are nonlinear. Moreover, the computational procedures of Bayesian inference may depend on the specific parametrization of dynamics models.

% Koopman
In this paper, we describe an approach toward uncertainty quantification that is parametrization-agnostic, using the \emph{transfer operator} description of dynamical systems. While commonly used in mathematical physics, the operator-theoretic view of dynamical systems has attracted attention for use in model order reduction, estimation, and control \cite{Mezic05}. Considering \emph{linear} operators (called transfer operators) over function spaces that represent the transition of observable or density functions, we can analyze, identify, and control nonlinear dynamical systems using linear techniques \cite{Mauroy15,Annoni16,Proctor16a,Surana16,Mauroy16,Mauroy16b,Surana17,Takeishi18}.

% Our contribution.
Model uncertainty quantification within the transfer operator view, however, is a relatively new research area, and developing suitable methods of transfer operator generation for robust optimization and control is the aim of our work.
To this end, we develop a method for sampling transfer operators using a metric \cite{ishikawa18} between nonlinear dynamical systems.
The proposed method is based on a Hamiltonian Monte Carlo in the space defined by the metric about a nominal model. Due to the transfer operator representation, our method requires only computing perturbations of linear operators (which are ultimately approximated as $d\times d$ matrices) to generate the uncertainty set, and thus does not require error propagation through complex parametrizations. Using both linear and nonlinear examples, we show how the sampling method preserves qualitative dynamical properties while efficiently searching dynamics space. We finally provide heuristics for sampling transfer operators with constraints. 

%==================================================

\section{RELATED WORK}

% model uncertainty
Uncertainty is inevitable as data are always finite and may include observation noise, but the intersection of uncertainty and the operator-theoretic view on dynamical systems has not yet been explored well.
Takeishi \emph{et al.} \cite{Takeishi17ijcai} discussed a probabilistic interpretation of a technique called dynamic mode decomposition (DMD), which is a widely used algorithm for computing transfer operators, but they only considered the uncertainty of spectral components of the operators.
Morton \emph{et al.} \cite{morton_deep_2019} considered the uncertainty of linear transition operators via the uncertainty of observable function values for a model based on neural networks.

In terms of parametric models, uncertainty and inference of parameters in nonlinear differential equations is a well-explored problem \cite{girolami2008bayesian}, with recent advances for accelerating Bayesian inference \cite{calderhead2009accelerating, niu2016fast}.  These methods are useful in certain contexts (e.g. where only parameters or spectral components are uncertain) but may not fully capture dynamics present in uncertain data. In this case, a more general method is needed to explore the space of dynamics, and this is where data-driven methods for approximating transfer operators such as extended DMD \cite{Williams15a} and observable eigenfunction learning \cite{korda2018optimal} are promising.

% state uncertainty
Complementarily to model uncertainty, state uncertainty has been at the center of interests in the control community.
Several researchers have investigated the computation of state uncertainty using operator-theoretic techniques  \cite{surana_linear_2016,schnitzer171109798,meyers19}. For example, state-space density propagation using Perron--Frobenius operators has been used for controller selection \cite{meyers19}. 

Our method is complementary to existing techniques for state- and parameter-space uncertainty quantification in two ways. First, rather than performing error propagation from ensembles of trajectories, we utilize a distribution induced by a metric between transfer operators -- thus it is a way to take an arbitrary set of initial models and expand it. Secondly, while the uncertainty set we produce can be used to extrapolate uncertain trajectories and perform state estimation, its use is not only limited to trajectory sampling but a wide range of optimization problems where a set of models is needed, without assuming a specific parametrization.

\section{BACKGROUND}

% In this section we review technical background on transfer operators and the pseudo-metric defined on them.

\subsection{Transfer operator theory of nonlinear dynamics}

A \textit{transfer operator} is a linear map $\mathcal{L}$ acting on functions $\phi : \mathcal{M} \to \C$ on phase space $\mathcal{M}$ of a dynamical system $F: \mathcal{M} \to \mathcal{M}$.
$\mathcal{L}$ is a linear operator over the function space $\{\phi\}$ and thus is a convenient alternative to the nonlinear state-space function $F$. 
% Importantly, it facilitates the study of $F$ via the operator $\mathcal{L}$:
% \[
% \begin{tikzcd}
%     \Phi \arrow{r}{\mathcal{L}} & \Phi \\
%     \mathcal{M} \arrow{u}{\phi} \arrow{r}{F} & \mathcal{M} \arrow{u}{\phi} 
% \end{tikzcd}
% \]

In computational fluid dynamics, statistical physics, control theory, and many other fields, the use of transfer operators in describing dynamical systems has seen a recent surge of popularity since it facilitates the usage of linear techniques (e.g. spectral methods) in analyzing nonlinear systems \cite{Mezic05}. 

When $\phi$ are phase-space densities $p_\mathcal{M}(x)$, $\mathcal{L}$ is referred to as the Perron--Frobenius operator $\mathcal{P}$ that acts as the pushforward operator for the Markov process $p_\mathcal{M}^+(x) = \mathcal{P}p_\mathcal{M}(x)$. When $\phi$ are arbitrary \textit{observables} $\phi: \mathcal{M} \to \C$, $\mathcal{L}$ is known as the Koopman (or composition) operator $\mathcal{K}$ \cite{Koopman255}.
% , named for Koopman and von Neumann's \cite{Koopman255} pioneering study on the spectrum of transfer operators for characterizing dynamical systems.
Seminal work by Mezi\'c \cite{Mezic05} showed the utility of transfer operator theory in data-driven system identification and model reduction.
Using Koopman eigenfunctions to obtain linear descriptions of nonlinear systems facilitated the study of dynamics via the spectrum of the Koopman operator (e.g., characterizations of chaos \cite{arbabi2017thesis}):
\begin{align*}
    \Phi(F^t(X)) = \Lambda^t\Phi(X),
\end{align*}
where $\Phi$ and $\Lambda$ are the eigenfunctions and the eigenvalues of the Koopman operator, respectively.
This led to the least-squares solutions for $\mathcal{K}$ by Schmid \emph{et al.} \cite{Schmid10}, called DMD. 

DMD has since been extended with dictionaries of nonlinear basis functions (extended DMD \cite{Williams15a}), observables in reproducing kernel Hilbert spaces \cite{Kawahara16,Klus19}, neural networks (e.g., \cite{Takeishi17neurips}), and the incorporation of control \cite{Proctor16a}.
In this paper, we build upon this prior work, exploring uncertainty quantification in the context of transfer operators which are learned from observation data. 

\subsection{Transfer operators for stochastic processes}

To discuss transfer operators in the context of uncertain data, let us move to a fundamentally probabilistic setting where observations $X_t \in \mathcal{M}$ represent a stochastic process. While the system may not be intrinsically noisy, the probabilistic view facilitates representing uncertainty about future states.
We adopt the framework of stochastic Koopman operators as introduced by Klus \emph{et al.} \cite{Klus19} and Song \emph{et al.} \cite{Song09}, and provide a brief background below.

Let a dynamical system $F: \mathcal{M} \to \mathcal{M}$ have invariant probability measure $\mu$, compactly supported over a measurable subspace $\Omega \subseteq \mathcal{M}$.
Let the sequence $\{X_t|t\ge 0\}$ be a Markov process with transition density given by $p(y | x) = \Pr\{F(X) = y | X = x\}$.
Given an observable $g \in L^2(\mathcal{M}, \mu)$, the action of the \emph{stochastic Koopman operator} $\koop$ on $g$ is defined by the conditional mean:
\begin{equation}
\label{eq:stk}
\koop g(x) = \E[g(F(X))|X = x] = \int_\Omega p(y|x)g(y) \mathrm{d}\mu(y).
\end{equation}
Its right-adjoint, the Perron--Frobenius operator $\mathcal{P}$, directly maps marginal distributions as $p^{+}(x) = \mathcal{P}p(x) = \int_\Omega p(y|x) p(x) \mathrm{d}\mu(x)$.
% \begin{equation*}
% p^{+}(x) = \mathcal{P}p(x) = \int_\Omega p(y|x) p(x) \mathrm{d}\mu(x).
% \end{equation*}

Let observables $\phi$ lie in an inner product space $\mathcal{H}$ ($\mathcal{H}$ is taken to be an RKHS in \cite{Klus19}, but for our purposes it can be the completion of a finite observable basis $|\Phi| = d$, and for trajectory visualization one may define these functions such that the preimage can be computed).
Now, we may define the Gramian $\mathcal{C}_{XX}$ (i.e., a \emph{cross-covariance operator} \cite{Song09}):
\begin{equation*}
\mathcal{C}_{XY} := \E_{XY}[\Phi(X)\otimes \Phi(Y)] = \int_{\Omega \otimes \Omega} \Phi(x) \otimes \Phi(y) \mathrm{d}p(x, y).
\end{equation*}

Using the relation $\mathcal{C}_{YX}f = \E_{Y|X}[f(Y)|X]\mathcal{C}_{XX}$ \cite{Song09}, we can express the stochastic Koopman and Perron--Frobenius operators \eqref{eq:stk} in terms of the cross-covariance operators as:
\begin{align}
\koop &\simeq C_{XY} (C_{XX} + \epsilon \mathcal{I})^{-1}, \\ 
\mathcal{P} &\simeq C_{YX} (C_{XX} + \epsilon \mathcal{I})^{-1},
\end{align}
that push forward observables $g\in \mathcal{H}$ and densities $p(x)$ on $\Omega$, respectively. $\epsilon$ is for ensuring invertibility. 

Considering perturbations to these operators, which are estimated from cross-covariances \emph{matrices} using finite trajectories, forms the population of our uncertainty set.

\subsection{Kernels over dynamical systems}
\label{sec:kernels}

% \TODO{mention binet-cauchy kernels?}
% \TODO{mention operator-valued kernels?}

To compare the behavior of two dynamical systems, various metrics exist in the literature. For Markov processes, total-variation distance between density functions is commonly used, and more general classes of linear metrics over Markov chains have been proposed \cite{DacaHKP16}. The standard operator norm can be used for bounded linear operators; when approximated by matrices, one can use an induced matrix norm, such as $\Tr(A^TB)$ where $A$ and $B$ are linear finite dynamical systems. However, both of these standard metrics fail to capture the iterated behavior of the respective systems. To address this, Viswanathan \emph{et al.} leverage the Binet--Cauchy theorem to define a kernel between trajectories or iterated maps \cite{vishwanathan2004binet}. 

In \cite{ishikawa18}, Ishikawa \emph{et al.} generalize the Binet-Cauchy kernel to general nonlinear dynamical systems defined by their respective Perron--Frobenius operators as follows.
For two dynamical systems $(D_1, D_2)$ specified by their initial values and maps $((X_{1,0}, f_1), (X_{2,0}, f_2))$, 
\begin{equation}\begin{aligned}
\label{eq:origkernel}
k_\mathcal{P}^{m, T} & ((X_{1,0}, f_1), (X_{2,0}, f_2)) \\ &:= \Tr\left( \bigwedge^m \sum_{t=0}^{T-1} (L_{1,h} \mathcal{P}_1^t \mathcal{I}_1)^*L_{2,h}\mathcal{P}_2^t \mathcal{I}_2\right),
\end{aligned}\end{equation}
where $\mathcal{I} : \C^n \to \mathcal{H}$ is an \emph{initial value operator} embedding initial data into a Hilbert space, $\mathcal{P}^t$ is the $t$th iterate of the Perron--Frobenius operator on $\mathcal{H}$, and $L_h : \mathcal{H} \to \mathcal{H}_\text{ob}$ embeds states into an observable space. There is a relation between \eqref{eq:origkernel} and the determinant kernel in \cite{vishwanathan2004binet}; please refer to \cite{ishikawa18} for details. In the proposed sampling method, we make use of this kernel \eqref{eq:origkernel} to generate perturbations to Koopman operators.

%==================================================

\section{PROPOSED METHOD}

In this section, we first define a kernel and a pseudo-metric over Koopman operators utilizing the previous studies on operator-theoretic kernels over dynamical systems \cite{Fujii17,ishikawa18}.
Then, we present sampling procedures for dynamical systems using the kernel, which depend on the Hamiltonian Monte Carlo method \cite{neal11hmc} with some heuristic modifications.

% \TODO{theory: suggestion for making the potential well-defined using discounting factor}
% \TODO{heuristics: initial conditional sampling from prior (uniform in our case)}
% \TODO{heuristics: spectral boundary reflection in leapfrog}

\subsection{Defining a kernel over Koopman operators}

Whereas the kernel proposed by Ishikawa \emph{et al.} \cite{ishikawa18} is defined for Perron--Frobenius operators on RKHS, we adapt it for the Koopman operator as follows.
Let the initial value and observable operators $\mathcal{I}, \mathcal{L}_h$ be already applied, and the kernel \eqref{eq:origkernel} evaluated on a dynamical system $\mathcal{K}$ which acts in this (observable) space.
Then we simplify \eqref{eq:origkernel} to:
\begin{align}
\label{eq:kernel}
k_\mathcal{K}^{m, T} & (f_1, f_2) := \Tr\left( \bigwedge^m \sum_{t=0}^{T-1} (\mathcal{K}_1^t)^* \mathcal{K}_2^t\right).
\end{align}
As noted in \cite{ishikawa18}, this kernel is convergent in the limit of $T$ for semi-stable $\mathcal{K}$, that is, those with spectral radius $\rho(\mathcal{K}) \le 1$ only.
To use $k_\mathcal{K}^{m, T}$ for general Koopman operators $\koop$, we use an exponential discounting scheme with factor $\lambda \geq 0$:
\begin{align}
\label{eq:discountedkernel}
    k_{\koop}^{m, T, \lambda} & (\koop_1, \koop_2) = \Tr\left( \bigwedge^m \sum_{t=0}^{T-1} e^{-\lambda t}(\koop_1^t)^*\koop_2^t \right).
\end{align}
Note that the discounting factor has been adopted also in previous studies on dynamical system kernels \cite{Vishwanathan07,Fujii17}.

Let us show the convergence of \eqref{eq:discountedkernel} as $T\to\infty$ informally, for a finite-dimensional approximation $K$ of $\koop$.
We first observe that $\sum_{t=0}^\infty A^t$ converges if $\lim_{t\to\infty}||A^t||_F = 0$ for any matrix $A$.
As the product of two convergent series is convergent, it suffices to show that $\lim_{t\to\infty}||e^{-\lambda t/2}K^t||_F = 0$ for all $K$.
Using Gelfand's formula,
\begin{align*}
    ||e^{-\lambda t/2}K^t||_F = e^{-\lambda t/2}||K^t||_F < e^{-\lambda t/2}\rho(K)^t.
\end{align*}
Even if $\rho(K) > 1$, the series is convergent if $\lambda > 2\log \rho(K)$.
Thus $k^{m, T, \lambda}$ converges for all $K$ and appropriate $\lambda$. 

\subsection{Sampling from distributions over transfer operators}
We propose a sampling procedure for dynamical systems models given an inner product defined over transfer operators.
We define a (pseudo-)metric bounded in $[0,1]$ for operators $K_1, K_2$, using a cosine similarity, as
\begin{align*}
    d_k(K_1, K_2) := \sqrt{1 - \frac{\langle K_1, K_2 \rangle_k^2}{\langle K_1, K_1 \rangle_k \langle K_2, K_2 \rangle_k}},
\end{align*}
where $\langle \cdot, \cdot\rangle_k$ denotes the inner product induced by a positive-definite kernel $k$.
As a baseline, we may also consider the standard linear kernel $k(A, B) = \Tr(A^\T B)$, which is not necessarily a good option for transfer operators.

Let $D \in [0, 1]$ be a random variable with density $p_D$.
Furthermore, let $K_0$ be a \textit{nominal} transfer operator (such as ones estimated by DMD).
Then we define a \textit{likelihood} of any dynamical system $K$ as:
\begin{align}
\label{eq:likelihood}
    \mathcal{L}(K \mid K_0) := p_D(d_k(K, K_0))
\end{align}
For example, we may assume $p_D$ is a Beta distribution with $\beta \gg \alpha$, or an exponential distribution with vanishing density past $1$.
Using this we may easily construct an uncertainty set of radius $r$ as $\Delta = \{ K \sim \mathcal{L}(K | K_0)\ |\ d_k(K, K_0) \le r \}$.

Sampling $K$ using $\mathcal{L}$ can be done in a number of ways. The conceptually simplest algorithm is rejection sampling: for any uniformly perturbed $K$, accept if $u < c\mathcal{L}(K | K_0)$ for $u \sim [0, 1]$ and convergence parameter $c$.
Unfortunately, even generating the initial uniform perturbations may be computationally intractable due to the high dimensionality of the samples.
This results in low acceptance rates for rejection sampling as well as random-walk MCMC methods such as Metropolis--Hastings.
We turn our focus to gradient-based MCMC methods which are able to generate distant proposals and achieve dimensionality-independent acceptance rates.

\subsection{High-dimensional sampling via HMC}
In his seminal work \cite{neal11hmc}, Neal introduced Hamiltonian Monte Carlo (HMC), which uses the gradient of the likelihood to simulate stochastic Hamiltonian dynamics whose stationary distribution is the posterior \eqref{eq:likelihood}.
It is well suited to our case since many kernels over dynamical systems are basically differentiable.

We adapt HMC for transfer operator sampling by introducing an auxiliary \textit{momentum} variable $R$ which is of the same dimension as $K$, whose relation to $K$ is given via Hamiltonian dynamics.
Let us define the \textit{potential} and the \textit{Hamiltonian} of an operator $K$ as:
\begin{align}
\label{eq:hamiltonian}
    U(K) &= -\log \Big[ p_D \big( d_k(K, K_0) \big) \Big] \\
    H(K) &= U(K) + \frac12 \Tr(R^TR)
\end{align}
Noting that the potential $U$ must be defined (and continuous) everywhere in order to achieve the correct stationary distribution, $d_k$ must be similarly well-behaved for all $K$.
Using the discounted kernel \eqref{eq:discountedkernel}, we generate samples $\{K\}$ of Koopman operators about a nominal $K_0$ via Hamiltonian dynamics in the potential defined by \eqref{eq:hamiltonian} using the \textit{leapfrog} integrator (we refer the reader to \cite{neal11hmc} for details).

\subsection{Heuristic: parallel HMC with uniform prior}
The mixing time of HMC is highly sensitive to the choice of discretization parameters (in particular, \texttt{n\_leapfrog}, \texttt{step\_size}).
In practice, Neal recommends $\epsilon \sim O(d^{1/4})$ \cite{neal11hmc}, and there also exist adaptive step-determining algorithms such as the No-U-Turn Sampler \cite{hoffman_2014_jmlr}; however, we find a tradeoff between computational efficiency (samples/step) and sufficient exploration of the model space when using HMC for transfer operators.
To compensate, we use a pre-run of HMC in a zero-potential to generate a uniform prior of samples, which are then used as initial conditions for parallel HMC sampling from the distribution \eqref{eq:likelihood}. We find that this accelerates the mixing time and wall-clock time for sampling significantly.

\subsection{Heuristic: HMC with spectral constraints}
In the interest of producing meaningful samples, one may wish to constrain the space of models using prior knowledge of the underlying system.
Constraints can be readily expressed as \textit{boundary conditions} on HMC without changing its stationary solution or reversibility (under some assumptions) -- known as Reflective HMC \cite{NIPS2015_5801}. 
As an example, suppose that we have knowledge that the system is structurally stable under any realistic perturbation; then, we can (roughly) encode this as a constraint on the spectral radius $\rho(K_0)-\epsilon \le \rho(K) \le \rho(K_0) + \epsilon$. 

Extending the reflective HMC procedure from \cite{NIPS2015_5801}, we describe a leapfrog integrator (Algorithm \ref{algo:reflectivehmc}) which ensures $f(K) \in [a, b]$ for any differentiable, scalar-valued $f$. In numerical experiments, we use $f(K) = \rho(K)$ and $[a, b] = [\rho(K_0) - .01, \, \rho(K_0) + .01]$.

\begin{algorithm}[b]
\caption{HMC with scalar constraints}
\label{algo:reflectivehmc}
{\small\begin{algorithmic}
\Procedure{BoundedLeapfrog}{$K,R,\epsilon, f, a, b$}
    \State $R\gets R - 1/2\epsilon\nabla_K U$
    \For{$i \gets 1$ to $L$}
        \While{$|\epsilon - \epsilon'| > \delta$} \Comment{Reflect till exhaustion}
            \State $\epsilon' \gets \epsilon$ 
            \State $K, R, \epsilon \gets \Call{Step}{{K},{R}, {\epsilon}, {f}, {a}, {b}}$
        \EndWhile
    \EndFor
    \State $R\gets R - 1/2\epsilon\nabla_K U$
    \State \textbf{return} $K, R$
\EndProcedure
\Function{Step}{$K,R,\epsilon, f, a, b$}
    \If{$f(K + \epsilon R) > b$}
        \State $K, \epsilon \gets \Call{FindMax}{{K}, {f}, {b}}$ \Comment{Find boundary}
        \State $R_\perp \gets \frac{\langle R, \nabla_K f\rangle}{\langle \nabla_K f, \nabla_K f\rangle} * \nabla_K f$ \Comment{Reflection plane}
        \State $R \gets R - 2R_\perp$
        \State \textbf{return} $K, R, \epsilon$
    \ElsIf{$f(K + \epsilon R) < a$}
        \State $\cdots$ \Comment{Defined similarly}
    \Else
        \State \textbf{return} $K + \epsilon R, R, 0$ \Comment{Otherwise, take full step}
    \EndIf 
\EndFunction
\end{algorithmic}}%
\end{algorithm}

\subsection{Computation of samples in practical settings}
We will give two formulations of HMC \eqref{eq:hamiltonian} for finite arguments, one explicitly over transfer operators, and one implicitly over observed trajectories. Either may be used.

\paragraph{Formulation 1}
Assume a transfer operator $\mathcal{K}$ for a dynamical system is approximated as a matrix $K \in \R^{d\times d}$.
Then, \eqref{eq:discountedkernel} simplifies to:
\begin{equation}
\label{eq:finitekernel}
\begin{aligned}
&k^{m, T, \lambda}_K(K, K_0)
\\&\quad = \sum_{I \subset [d], |I| = m} \det \left(\sum_{t=0}^{T-1} e^{-\lambda t}(K^t)^\T K_0^t\right)_{[I, I]},
\end{aligned}
\end{equation}
where $A_{[I, I]}$ denotes the submatrix given by indices $I$. 
This formulation can be used when we have a nominal Koopman operator estimation $K_0$ and want to generate perturbed $K$'s.

\paragraph{Formulation 2}
Assume observations $x_i \in \R^d$ of a dynamical system are given as a matrix $X$.
Then, we may define the Hamiltonian \eqref{eq:hamiltonian} for trajectories as:
\begin{align*}
    H(X) &= -\log\left[p_D(d_k(X, X_0))\right]  + 1/2\langle R, R\rangle.
\end{align*}
The kernel \eqref{eq:discountedkernel} can be defined for trajectories of length $T$ in a similar fashion to \eqref{eq:finitekernel}.
With $T$ fixed, an explicit discounting term is no longer needed.
As an example, for $m=2$: 
\begin{align*}
    &k^{2, T, \lambda}(X, X_0) = \sum_{i,j\in [1, T]} \det \begin{bmatrix} k(x_{i},x_{j}) & k(x_{i},x_{0,j}) \\ k(x_{0,i},x_{j}) & k(x_{0,i},x_{0,j})\end{bmatrix}
\end{align*}
for some \textit{feature} kernel $k(\cdot,\cdot)$.
Here, $x_i$ and $x_{0,i}$ are the elements of $X$ and $X_0$, respectively.
If we have explicit observables $\phi$, simply let $k(x, y) = \langle\phi(x), \phi(y)\rangle$.
This formulation allows one to sample trajectories $\{X\}$ about some observed trajectory $X_0$. When used in conjunction with a parametrized dynamics model $F_\theta$, \{X\} can be used for inference procedures (e.g., E-M) on $\theta$.

\begin{figure*}[t]
    \vspace*{1mm}
    \centering
    \includegraphics[width=0.98\linewidth,height=6cm]{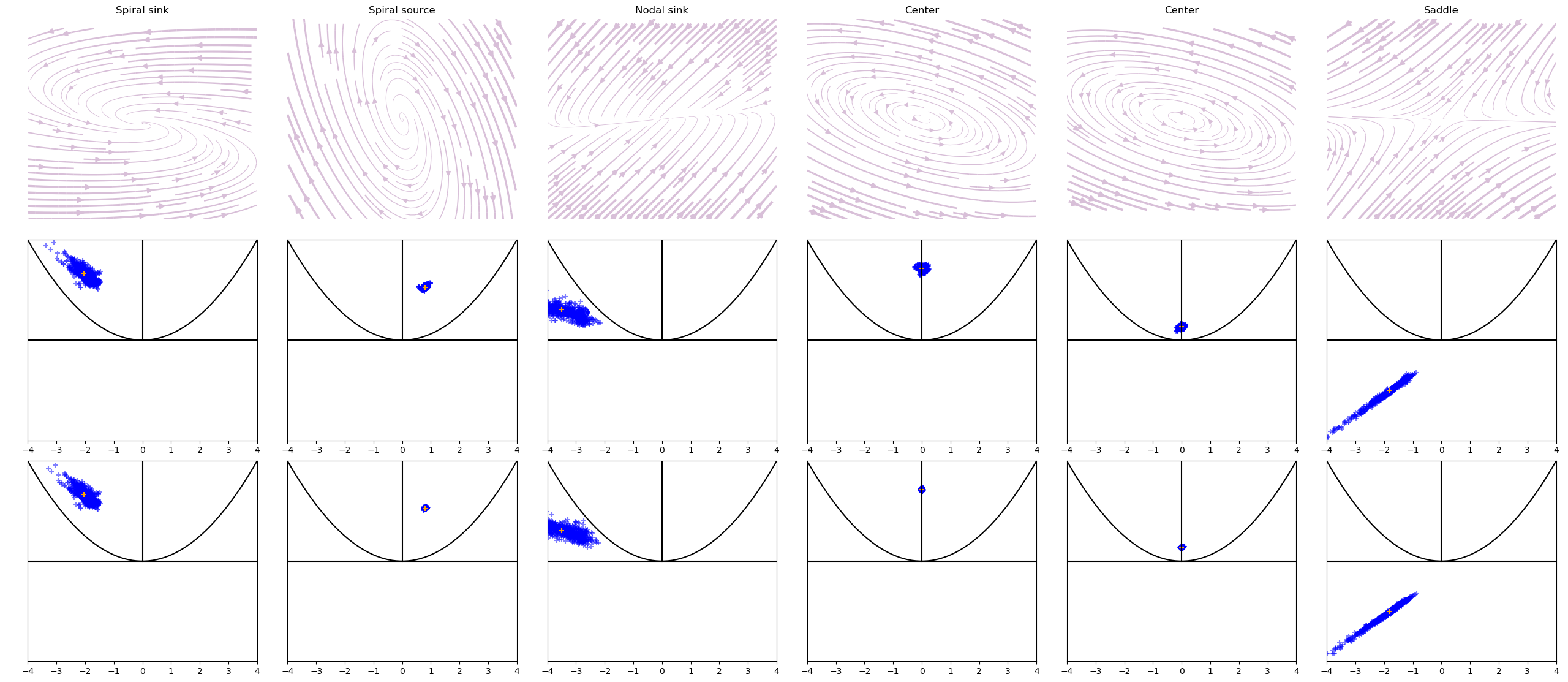}
    \caption{Perturbations of 2x2 systems. Left to right: spiral sink, spiral source, nodal sink, center, center, saddle. Top to bottom: trace-determinant plots for baseline vs. proposed method ($m=2, T=80, \lambda=2\log\rho(A_d)$).}
    \label{fig:2x2}
\end{figure*}

\section{NUMERICAL EXAMPLES}

We show how our sampling procedure compares against baseline perturbation methods in generating meaningful perturbations of dynamical systems\footnote{Codes: \url{github.com/ooblahman/koopman-robust-control}.}.
For example, methods in robust optimization (RO) take the general form 
\begin{align*}
    \min_\theta \max_{A \in \Delta} J(A, \theta) \quad\text{where}\quad \Delta := \{A + \Delta A\},
\end{align*}
where $A$ is a model and $\Delta A$ is an \textit{uncertainty structure}.
$\Delta A$ may have a particular form, e.g. block-diagonal, or unstructured, e.g. $\Delta A \sim \mathcal{N}(0_{d\times d}, \Sigma^2)$.
The goodness of the RO minimizer depends solely on the choice of uncertainty set.
Either of these uncertainty structures essentially induces a distribution over the \textit{norms} of perturbations; we note that this is an assumption, and the subject of our testing is whether this assumption is valid when it is known the perturbed matrices represent dynamical systems.

A norm-bounded perturbation set essentially translates as sampling with a trace kernel:
\begin{align}\label{eq:tracekernel}
    k(A, A_0) = \Tr(A^\T A_0),
\end{align}
which we use as a baseline for comparison. As mentioned in section \ref{sec:kernels}, despite its naturalness, it is not necessarily suitable for dynamical systems.
By contrast, the proposed uncertainty structure induced by \eqref{eq:finitekernel} incorporates the iterated behavior of $K$. 
We demonstrate that the latter better preserves key properties of dynamical systems such as structural stability and attractor basins, while effectively exploring dynamics space, on both linear systems of ODEs and nonlinear systems via the Koopman operator. 

\subsection{2-dimensional LTI systems}
We first consider linear systems of the form $\dot{x} = Ax$ via discretization as $A_d = e^{A \Delta t}$.
We use simple 2x2 systems in order to clearly characterize the dynamics in a trace-determinant plot.
Using the discounted kernel \eqref{eq:finitekernel}, we are able to generate perturbations of both source- and saddle-types in addition to the semistable regimes.
We use spread parameter $\beta = 5$, HMC step $\epsilon = 10^{-4}$, HMC leapfrog $L = 100$, and generate $N = 1000$ samples with $k = 50$ initial conditions for every test shown.
We compare the following two kernels:
\begin{align*}
    \Tr(A^\T A_0)~\text{in}~\eqref{eq:tracekernel} \quad\text{and}\quad k^{m, T, \lambda}_K(A, A_0)~\text{in}~\eqref{eq:finitekernel},
\end{align*}
termed the ``trace kernel'' and ``Koopman kernel'', respectively.

The strength of the proposed method, using $k^{m, T, \lambda}_K$, can be seen when the nominal system is within a region of structural instability (see the two \textit{center} systems, Figure~\ref{fig:2x2}).
The trace kernel perturbations venture easily into spiral sink or spiral sources, which are distant in dynamics terms but very close in absolute norm.
In these cases, the proposed method using the Koopman kernel retains a much tighter spread in the trace-determinant plane.
Furthermore, it can be seen from the posterior distributions that the proposed method is able to explore distant dynamics while staying bounded within structurally similar regions.

%\FloatBarrier
\subsection{Unforced Duffing oscillator}

\begin{wrapfigure}[7]{r}{3cm}
    \centering
    \vspace*{-10mm}
    \includegraphics[clip,width=2.8cm]{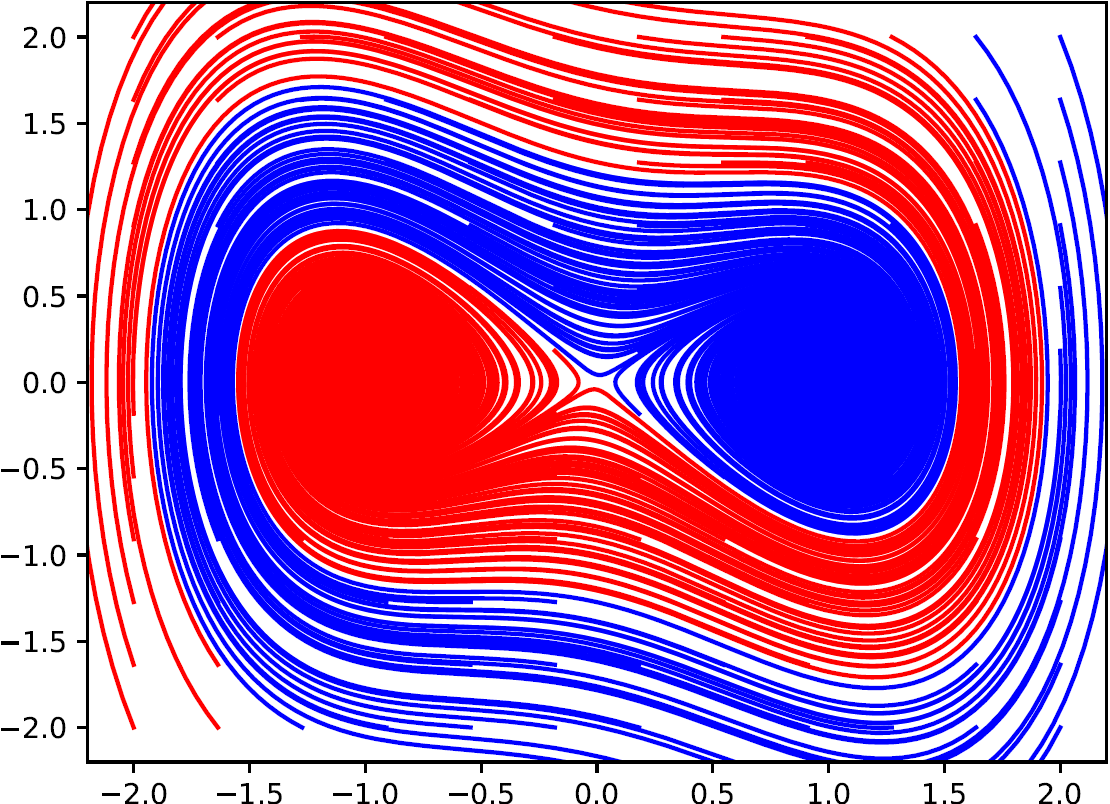}
    \vspace*{-1mm}
    \caption{Basins of attraction of \eqref{eq:duffing}.}
    \label{fig:duffing}
\end{wrapfigure}
Next, we consider perturbations of a nonlinear system.
We use the unforced Duffing equation:
\begin{equation}
    \label{eq:duffing}
    \dot{x} = y, \quad \dot{y} = -0.3y + x - x^3,
\end{equation}
whose basins of attraction are like in Figure~\ref{fig:duffing}.
\begin{figure*}
    \centering
    % \subfloat[]{
        \includegraphics[width=0.48\linewidth]{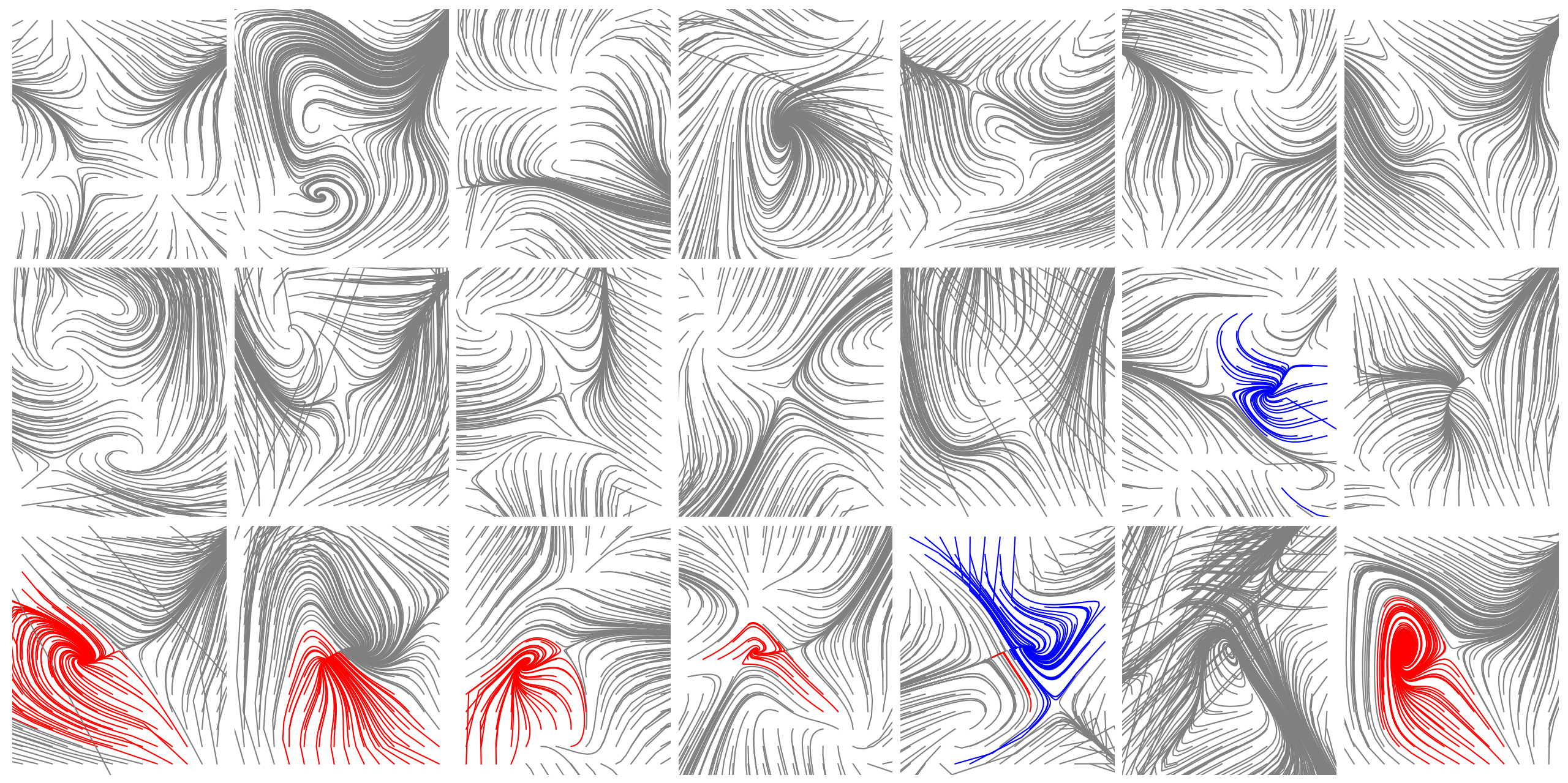}
        % \label{fig:duffing_baseline}
    % }
    \hspace{1em}
    % \subfloat[]{
        \includegraphics[width=0.48\linewidth]{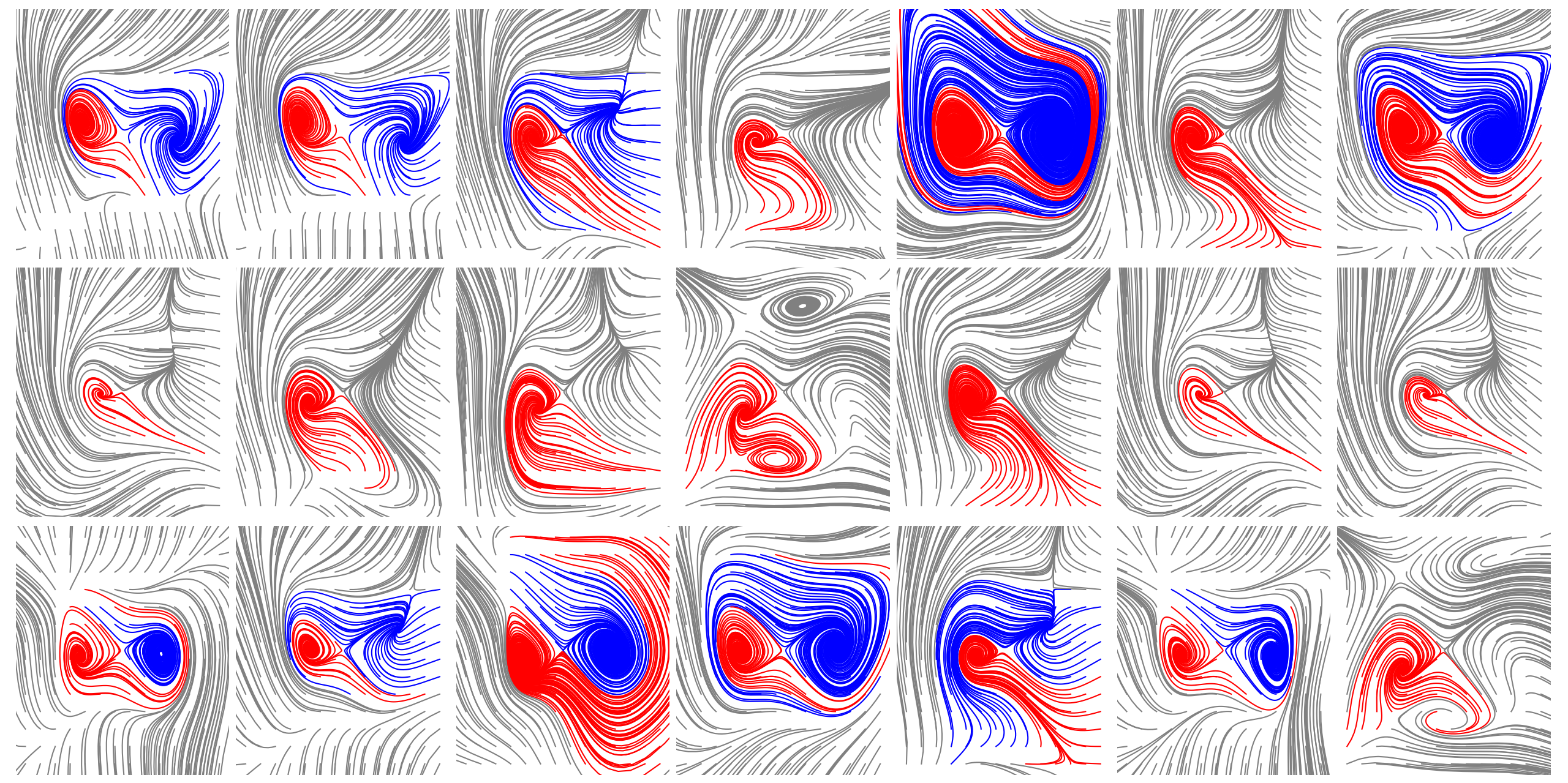}
        % \label{fig:duffing_kernel}
    % }
    \caption{(\emph{left}) Perturbations via the trace kernel and (\emph{right}) perturbations via the Koopman kernel ($T=80, m=2, \lambda=0$).}
    \label{fig:duffing_result}
\end{figure*}
We use simulated trajectories of length $t = 400$ seconds with $8000$ samples per trajectory across $144$ initial conditions in the range $[-2, -2], [2, 2]$.
Using $15$ polynomial observables with maximum degree $5$, we compute the Koopman operator $K$ as $K = YX^\T (XX^\T)^\dagger$.
For all experiments we use spread $\beta = 5$, HMC step $\epsilon = 5\times10^{-5}$, HMC leapfrog $L = 200$, $N = 2000$ samples, and $k = 200$ initial conditions.
In the shown perturbations, trajectories in the left and right basins of attraction are highlighted in red and blue, respectively (Figure~\ref{fig:duffing_result}).

We immediately observe some qualitative differences between the two perturbation sets.
First, it is apparent that perturbations via the Koopman kernel preserve attractor structure in most samples, versus almost none in the baseline setting. In the case of the Duffing oscillator, this is a defining feature, and such preservations are important consideration for any robust prediction or control procedure over dynamics models.
Second, a large proportion of samples produced by the baseline method are diverging; these would need to be manually filtered out if used in a robust optimization setting.
We observe in experiments that this can be mitigated by restricting the norm of perturbations (i.e., increasing $\beta$), but this comes at the cost of decreased exploration of dynamics space and changes the robustness of an RO solution using the perturbation set (moreover, manual filtering changes the posterior, altering the RO problem).
We also find that the spectral radius constraint alleviates many of these concerns with the baseline method, however, non-convex reflection is not a trivial procedure to implement in HMC and has not been a typically used method in generating uncertainty sets.

Finally, while attractors are mostly preserved in our method, the attractor basins seem to undergo some geometry warping. This suggests an interpretation of our perturbation method as warping the underlying potential wells, which may have meaningful physical interpretations. 

%\subsection{Robust control of a Duffing oscillator}

%%%%%%%%%%%%%%%%%%%%%%%
% \TODO{prediction with interval (maybe with low-dimensional dynamics), comparing to EDMD/KDMD (without uncertainty), Bayesian linear dynamical system (which is basically equivalent to Bayes DMD), and GP dynamics model? Data: FitzHugh-Nagumo, Van der Pol, Lorenz, and some real-world data (CMU MOCAP?)}
% http://mocap.cs.cmu.edu/
% \TODO{control w/ Burger's Equation, maybe}
%
% 1 ... proof-of-concept experiment ... generate (simulate) trajectories from nominal and perturbed operators ... maybe look at eigenvalues of them ...
% baseline = perturbation by Euclidean distance on matrix elements, perturbation on trajectories (and then matrix estimation), and maybe GP (base on e-DMD; try low-dimensional feature space containing full state observable)
%
% (maybe: examine baseline uncertainty set near bifurcation for simple system, vs. our method), e.g. any non-chaotic system with bifurcations (pendulum)
%
% 2 ... prediction with interval ...
% baseline=?...
%
% 3 ... robust MPC ... compare with eDMD-based MPC?, maybe time / amount of disturbance / H-inf norm ...
%%%%%%%%%%%%%%%%%%%%%%%

%==================================================

\section{CONCLUSIONS}

In this work, we developed a method for sampling from distributions over dynamical systems using transfer-operator representations, leveraging operator-theoretic metrics to generate perturbations.
We suggested to use the method for model uncertainty set generation, which is a universal problem in robust control and prediction and an important step for uncertainty quantification.
Future directions of research include expressing constraints over sampled dynamical systems where we may have domain-specific knowledge.

% - manual methods of constructing uncertainty sets is hard. naive perturbations of matrices easily result in dynamics that don't make sense
% - directly sampling distributions over transfer operators, with meaningful constraints, is possible with MCMC methods. automatic differentiation packages make this particularly easy
% - probabilistic methods for estimating dynamical models (monte carlo system identification) allow the quantification of and optimization for uncertainty in the case of noise, incomplete data, or incorrect models
% - future methods, we would like to take the idea of expressing constraints over sampled dynamical systems where we may have domain-specific knowledge. For example, in the case of Duffing equation, we keep the potential well bounded. In the case of van der pol oscillator, we retain the property of limit cycle. In other dynamical systems such as neural nets, attractor nets, etc, we may have interesting constraints e.g. weight sharing type schemes, or ways to prevent stuff like gradient exlposion in more principled way. Transfer learning?

%==================================================

%\addtolength{\textheight}{-12cm}

%==================================================

\bibliography{refs,tmp}
\bibliographystyle{ieeetr}

% \section*{ACKNOWLEDGMENT}
%
% ...

%==================================================

%\bibliographystyle{IEEEtran}
%\bibliography{}

\end{document}